\documentclass[11pt,a4paper]{amsart}

\usepackage{latexsym}

\usepackage[dvips]{graphics}
\usepackage{amssymb}
\usepackage{amsthm}
\usepackage{amsmath}
\usepackage{amstext}
\usepackage{amscd}
\usepackage{color}

\usepackage{ifmtarg}
\makeatletter
  \newcommand{\TODO}[1]{\@ifmtarg{#1}{\emph{\textcolor{red}{\textbf{TODO}}}~}{\textcolor{red}{ \emph{\textbf{TODO:}~#1~}}}}
  \newcommand{\FORUS}[1]{\@ifmtarg{#1}{\emph{\textcolor{blue}
  {\textbf{TODO}}}~}{\textcolor{blue}{ \emph{\textbf{FOR US:}~#1~}}}}
\makeatother

\numberwithin{equation}{section}

\theoremstyle{plain}
\newtheorem{thm}{Theorem}[section]
\newtheorem{lem}[thm]{Lemma}

\newtheorem{theorem*}{Theorem}[]

\newtheorem{lemma}[thm]{Lemma}
\newtheorem{prop}[thm]{Proposition}
\newtheorem{cor}[thm]{Corollary}

\theoremstyle{definition}
\newtheorem{defi}[thm]{Definition}

\newtheorem{rmk}[thm]{Remark}

\theoremstyle{definition}
\newtheorem{defn}[thm]{Definition}

\theoremstyle{remark}
\newtheorem{rem}[thm]{Remark}

\newtheorem{question}{Question}

\newcommand{\C}{\mathbb{C}}
\newcommand{\N}{\mathbb{N}}

\newcommand{\h}{{h}}

\newcommand{\kk}{{k}}

\newcommand{\D}{D}
\newcommand{\DD}{\Delta}

\DeclareMathOperator{\id}{id}

\newcommand{\impli}{\varphi}
\newcommand{\derimpli}{\psi}

\newcommand{\contour}{\mathcal C}
\newcommand{\polar}{\mathcal P}

\newcommand{\point}{P}

\DeclareMathOperator{\sing}{Sing \, }

\newcommand{\my}{\mathfrak m}

\usepackage{enumitem}
\setlist[itemize]{labelindent=.6em, itemindent=1em, leftmargin=!, label=\textbullet}
\usepackage{ifmtarg}

\newcommand{\snt}{$\nu$-transverse}
\newcommand{\unt}{$\nu$-transverse }

\newcommand{\PW}{\mathcal {PW}}

\newcommand{\dimtype}{\operatorname {d.t. }}

\title{Zariski's dimensionality type of singularities. Case of dimensionality type 2}

\author{Adam Parusi\'nski and Lauren\c tiu P\u aunescu}
\address {Universit\'e C\^ote d'Azur,  CNRS,  LJAD, UMR 7351, 06108 Nice, France}
\email{adam.parusinski@univ-cotedazur.fr}

\address{School of Mathematics and Statistics, The University of Sydney,
  Sydney, NSW, 2006, Australia }%
\email{laurentiu.paunescu@sydney.edu.au}%

\thanks{The first author is grateful for the support and hospitality of the Sydney Mathematical Research Institute (SMRI). 
Partially supported by ANR project LISA (ANR-17-CE40-0023-03).
}

\begin{document}
\begin{abstract}
In the 1970s O. Zariski  introduced a general theory of equisingularity 
for algebroid and  algebraic hypersurfaces over an algebraically closed field of characteristic zero.
His theory builds  up on understanding  the dimensionality type of hypersurface singularities, notion defined recursively by considering the discriminants loci of successive "generic" corank $1$ projections. 
The theory of singularities of dimensionality type 1, that is the ones appearing generically in codimension 1, was developed by Zariski in his foundational papers on equisingular families of plane curve singularities. 
In this paper we completely settle the case of dimensionality type 2,  by studying Zariski  equisingular families of surfaces singularities, not necessarily isolated, in the three-dimensional space.

\end{abstract}
\maketitle
\tableofcontents


\section{Introduction}

O. Zariski \cite{Zariski1979} introduced a general theory of equisingularity 
for algebroid and  algebraic hypersurfaces over an algebraically closed field of characteristic zero.  This theory is based on the notion of dimensionality type of hypersurface singularities introduced and developed by Zariski in loc cit..   The dimensionality type is defined recursively by considering the discriminants loci of successive "generic" corank $1$ projections. The singularities of dimensionality type 0 are the regular points of the hypersurface.  The singularities   of dimensionality type 1 are isomorphic to  the total space of an equisingular family of plane curve singularities, corresponding to the points of dimensionality type 0 of the discriminant of a generic projection. The singularities   of dimensionality type 2  correspond to the points of dimensionality type 1 of the discriminant of a generic projection, see definitions \ref{def:dimensionalitytype} and \ref{def:generic equisingularity} for a precise formulation. 
In general, the points of a fixed dimensionality type of an algebraic hypersurface form a smooth locally closed subvariety and thus the singularity type defines a canonical (independent of the system of coordinates) stratification of the hypersurface.  
For more on Zariski equisingularity including a historical account see \cite{Lip97}.  For a recent survey see \cite{MR4367438}.

The notion of a generic projection plays a crucial role in Zariski's approach but his very definition of such a projection involves adding all the coefficients, of a generic local formal change of  coordinates,  as indeterminates to the ground field.  This is not necessary in the case of  dimensionality type 1, the case when the transverse projections are sufficiently generic for Zariski's purpose.  On the other hand, this is no longer the case for the singularities of dimensional type 2.   Indeed, in  \cite{luengo85} Luengo gave an example of a family of surface singularities in $\C^3$ that is Zariski equisingular for one transverse system of coordinates but not for a generic linear system of coordinates.  As Zariski showed in  \cite{Zariski1979},  a generic polynomial projection is sufficient for any dimensionality type, though he gave no explicit bound on the degree of such polynomial map. This makes an algorithmic computation of Zariski's canonical stratification impossible. The algebraic case was studied in more detail by Hironaka \cite{hironaka79}, where the semicontinuity in Zariski topology of such a degree is shown.  

The question whether a generic linear projection 
is always sufficient is still open for dimensionality type at least three.  The case  of dimensionality type 1 follows from the above mentioned work of Zariski. 
An important progress in the case of dimensionality 2 was done in \cite{brianconhenry80}, where generic  Zariski equisingular families of complex analytic isolated surface singularities  in the 3D-space were considered.  It was shown in 
\cite{brianconhenry80} that, for such families,  the generic Zariski equisingularity is equivalent to the constancy of: Teissier's numbers (multiplicity, Milnor number, and Milnor number of a generic plane section), the number of double points and the number of cusps of the apparent contours of the generic projection of the generic  fibre of a mini-versal deformation.  All these numbers are local analytic diffeomorphism invariants, and therefore a linear generic change of coordinates is generic in the sense of Zariski.

We introduce in Section \ref{sec:main} the nested uniformly transverse (\snt) Zariski equisingularity that implies both generic and generic linear Zariski equisingularities. Now we state our main theorem.  Likewise  in \cite{Zariski1979},  this theorem is stated for the algebroid varieties, that is the ones defined over the ring of formal power series.

\begin{thm}\label{thm:mainequivalence}
Let $(V,\point)\subset (\kk^{r+1},\point )$ be an algebroid hypersurface and let $S$ be a nonsingular subspace of $(V,\point)$ of dimension $r-2$.   Assume that the dimensionality type of $V$ at $P$ is at least 2. Then the following conditions are equivalent:
\begin{enumerate}
  \item 
  $V$ is \unt Zariski equisingular along $S$ at 
$\point$.  
\item
For every $Q\in S$ the dimensionality type of $V$ at $Q$ 
is equal to $2$. 
\item
$V$ is generic Zariski equisingular along $S$ at $\point$. 
\item
For a local system of coordinates at $\point$,  $V$
 is generic linear Zariski equisingular along $S$ at $\point$. 
 \item
For all local systems of coordinates at $\point$,  $V$
 is generic linear Zariski equisingular along $S$ at $\point$. 
\end{enumerate}
\end{thm}

Theorem \ref{thm:mainequivalence} is proved in Section \ref{ssec:proofMainThm}. 
The equivalence (2)$\Longleftrightarrow $(3) already follows  from  \cite{Zariski1979}.  Our main results are 
(1)$\Longrightarrow$(3) and (1)$\Longrightarrow$(5),  whose  proofs are based on two main technical results, Theorem \ref{thm:genericwithparameter} and Proposition \ref{prop:retionalinparameter}.  The implication (4)$\Longrightarrow$(1)
is shown in Lemma \ref{lem:(4)==>(1)}.  Consequently we get two important implications (4)$\Longrightarrow$(3) and (4)$\Longrightarrow$(5) that we do not know how to prove directly.  The other implications to complete the proof of Theorem \ref{thm:mainequivalence} 
are standard, see Section \ref{ssec:proofMainThm}.  

Proposition \ref{prop:retionalinparameter} allows us to translate Zariski equisingularity with respect to a global parameter, 
for instance given by a finitely dimensional family of linear or polynomial automorphisms, to Zariski equisingularity with respect 
to a local parameter. Its statement is valid in any dimension.  

The notion of $\nu$-transversality can be generalized to higher dimensionality types but our proof of  Theorem \ref{thm:genericwithparameter},  that shows that local deformations of \unt Zariski equisingular families are themselves Zariski equisingular,  is, at the moment, restricted to the families of surface singularities.   
Its proof, see Section \ref{sec:proof}, is based on a precise analysis of the  \unt Zariski equisingular  
families of surface singularities in the three-dimensional space 
and on the parameterizations of the associated families of 
polar curves and their neighborhoods (polar wedges).  
We recall it in Section \ref{sec:parameterizations}.  This technique stems from \cite{brianconhenry80} and \cite{teissierRabida82} and was recently developed in \cite{BNP2014} and  \cite{NPpreprint}.  We applied it in \cite{PPpreprint19} to show that such families of surface singularities, in the complex analytic case, are bi-Lipschitz trivial  by constructing their explicit Lipschitz stratification.  


Our results are stated over the ring of formal power series.  Sometimes, following Zariski, we use the word analytic 
meaning defined by formal power series.  
We briefly discuss  the algebraic case and the complex analytic case in Sections \ref{ssec:algebraic} and \ref{ssec:analytic} 
respectively.

\subsection{Notation} 
By $\kk$ we denote an algebraically closed field  of characteristic zero that is fixed throughout this paper.

For a polynomial monic in $z$, 
$F(x,z) = z^d+ \sum_{i =1}^d a_i(x) z^{d-i}$, with coefficients analytic functions in $x$, we denote by $\D_F(x)$ its discriminant, and by 
$\DD_F$ its discriminant locus, the zero set  of  $\D_F$.


\section{Zariski equisingularity along a subspace and the dimensionality type}
\label{sec:dim-type}

Let $\kk$ be an algebraically closed field of characteristic zero.
Consider an algebroid  hypersurface $V=f^{-1}(0) \subset (\kk^{r+1},\point )$ at $\point \in \kk^{r+1}$ defined in a local system of coordinates by a formal power series $f\in \kk [[x_1, \ldots , x_{r+1}]]$ that we always assume  reduced.  Let 
$S$ be a nonsingular algebroid subspace of $\sing V$. In 
\cite[Definition 3] {zariski71open} Zariski introduced the notion of  algebro-geometric equisingularity of $V$ at $\point $ along $S$, now called Zariski equisingularity.  
It is defined recursively with respect to the codimension  of $S$ 
in $V$ by taking the images and the discriminants of some successive co-rank 1 projections.  

We say that a local projection $\pi: (\kk^{r+1},\point ) \to (\kk^{r},\pi(\point ))$ is \emph{permissible} if it is of rank $r$, the restriction $\pi_{|V}$ is finite, and the fiber $\pi^{-1}(\pi (\point ))$, 
that is a non-singular curve,  is not tangent to $S$. 
If $\dim S =\dim V - 1$ then we say that $V$ is 
\emph{Zariski equisingular at $\point $ along $S$} if there is a local permissible projection $\pi$ such that $\pi(\point )$ 
is a regular point of the discriminant locus $\Delta$  of $\pi_{|V}$.   In the general case,  
we say that $V$ is \emph{Zariski equisingular at $\point $ along $S$} if there is a local permissible projection $\pi$, such that the discriminant locus $\Delta$ of $\pi_{|V}$ is Zariski equisingular at $\pi(\point )$ along $\pi (S)$ (that implies that $\pi (S)$ is a nonsingular subspace of $\sing \Delta$) or  $\pi(\point )$ is a regular point of $\Delta$.  
 
  Stronger notions of Zariski equisingularity are obtained 
   if one assumes $\pi^{-1}(\pi (\point))$ is not tangent to $V$ at $\point$ 
   (\emph{transverse Zariski equisingularity}) or that $\pi$ is generic in the sense of \cite{Zariski1979}, that we will recall next. The 
   \emph{generic projection} of \cite{Zariski1979} is the map 
    $\pi_u (x) = (\pi_{u,1}(x), \ldots, \pi_{u,{r}}(x))$, with 
  \begin{align}\label{eq:genprojection}
\pi_{u,i} (x) = \sum_{d\ge 1} \sum_{\nu_1+ \cdots \nu_{r+1}=d} u^{(i)}_{\nu_1, \cdots, \nu_{r+1} } x^\nu ,  
  \end{align}
 considered as a projection defined over $\kk^*$, any field extension of $k$ that contains all coefficients $u^{(i)}_{\nu_1, \cdots, \nu_{r+1} }$ as indeterminates. 
The notions we are going to define now are independent of the choice of 
this field extension, see  \cite{Zariski1979} and \cite{zariski1980}. 

 Let $\Delta_u ^*\subset ({(\kk^*)}^{r},\point )$ denote the discriminant locus of 
 $\pi_u|_{V^*}$, where $V^*=f^{-1}(0) \subset ({(\kk^*)}^{(r+1)},\point )$, $\point_0^*= \pi_u (\point)$, $S_0^* = \pi_u (S)$.

\begin{defn}\label{def:generic equisingularity}
We say that $V$ is \emph{generic Zariski equisingular at $\point $ along $S$} 
 if $\Delta_u ^*$ is generic  Zariski equisingular 
at $\point_0^*$  along $S_0^*$  (or $\point_0^*$ is a regular point of $\Delta_u ^*$).
\end{defn}  
Let us recall from  \cite{Zariski1979} the related notions of 
dimensionality type and stratification by the dimensionality type. 
The definition of dimensionality type is again recursive. It is defined 
for any point $Q$ of $V$, not only for the closed point $\point $. 

\begin{defn}\label{def:dimensionalitytype}
Any simple point $Q$ of $V$ is  of dimensionality type $0$.  
Let $Q$ be a singular point of $V$ and let $Q^*_0= \pi_u (Q)$.  
Then \emph{the dimensionality type of $V$ at $Q$}, 
denoted by $\dimtype (V,Q)$, is equal to $1+ \dimtype (\Delta^*_u,Q^*_0)$.  
\end{defn}

The set of points where the dimensionality type is constant, say equal to $\sigma$, is either empty or a nonsingular locally closed subvariety of $V$ of codimension $\sigma$, \cite[Theorem 5.1]{Zariski1979}.   
A  singularity is of dimensionality type $1$ if and only if  it is isomorphic to the total space of an equisingular family of plane curve singularities, see \cite[Theorem 4.4]{zariski65-S2}. In general, the dimensionality type defines a canonical stratification of $V$ that satisfies the frontier condition. By  
\cite[Appendix]{Zariski1979} and by \cite{zariski75}  $V$ is equimultiple along each stratum.   
Then $V$ is generic Zariski equisingular along $S$ at $\point$ if and only if  $S$ is contained in the stratum containing $\point$.

In \cite[Proposition 5.3]{Zariski1979} Zariski shows that, with  the above definitions, the generic projection $\pi_u$ can be replaced by 
a condition that involves almost all projections 
$\pi_{\bar u}: \kk^{r+1}\to \kk^{r}$, as explained below.    
 One says that a property holds \emph{for almost all projections} if 
 there exists a finite set of polynomials $\mathcal G=\{G_\mu\}$ in the indeterminates  $u^{(i)}_{\nu_1, \cdots, \nu_{r+1} }$ and coefficients in $\kk$ such that this property holds for all projections $\pi_{\bar u}$ for $\bar u$ 
 satisfying  $\forall_\mu G_\mu (\bar u)\ne 0 $.  
 Here the bar denotes the specialization $u \to \bar u$, i.e. 
 we replace all indeterminates $u^{(i)}_{\nu_1, \cdots, \nu_{r+1} }$ by 
 elements of $\kk$, $\bar u^{(i)}_{\nu_1, \cdots, \nu_{r+1} }\in \kk$.  Thus, for almost all  projections $\pi_{\bar u} $ 
 the dimensionality type $\dimtype (V,\point )$ equals 
 $1+ \dimtype (\Delta_{\bar u},\pi_{\bar u} (\point ))$, where $\Delta_{\bar u}$ 
 denotes the discriminant locus of ${\pi_{\bar u}}_{|V}$.  
 Since the finite set of polynomials $\mathcal G$ involves nontrivially only finitely many indeterminates  $u^{(i)}_{\nu_1, \cdots, \nu_{r+1} }$, we may specialize the remaining ones to $0$, and then the projection 
 $\pi_{\bar u} $ becomes polynomial. 
 This means that, as soon as we know the set of polynomials $\mathcal G$, 
  we may compute the dimensionality type of $V$ at $\point $ 
 just by computing $\dimtype (\Delta_{\bar u},\pi_{\bar u} (\point ))$, 
 for only one polynomial projection $\pi_{\bar u} $, satisfying $\forall_\mu  G_\mu (\bar u)\ne 0 $.  
 Similarly, in order to check whether the generic Zariski  equisingularity of $V$ along $S$ holds at $\point $, 
 it suffices to check it for $\Delta_{\bar u} ^*$  along $\pi_{\bar u} (S)$ at $\pi_{\bar u} (\point )$. 
 In the sequel we shall omit the the bar notation $\bar u$ whenever it is clear from the context.

In general, no specific information on the polynomials of 
$\mathcal G$ is given in \cite{Zariski1979}.  In particular, no 
explicit bound on their degrees or on the weight $|\nu|=\sum \nu_k$ of the variables  
$u^{(i)}_{\nu_1, \cdots, \nu_{r+1} }$ they depend on. 
  Only the case of dimensionality type $1$ is known.  In this case, 
  by \cite{zariski65-S2}, the   
 generic and the transverse Zariski equisingularities coincide.  
 This means that one may choose all $G_\mu$ to depend only on the variables of weight 
$1$, $u^{(i)}_{\nu_1, \cdots, \nu_{r+1} }$ such that $|\nu|=1$, and   
such that $\mathcal G$ expresses the conditions that the projection 
$\pi_{\bar u}$ is linear of rank $r$ and the leading homogeneous form of $f$ does not vanish on its kernel.  
The above property is special for the dimensionality type $1$, the transversality of the projection does not imply 
 the genericity for the dimensionality type $2$.  
In \cite{luengo85} Luengo gave an example of a family of surface 
singularities in $\C^3$ that is Zariski equisingular for one transverse 
projection but not for  generic ones.  

\begin{rem}\label{rem:comment1}
It is clear, by its definition,  that the notion of generic Zariski equisingularity does not depend on the choice of a local system of coordinates $x=(x_1, \ldots, x_{r+1})$.  One can also see it 
directly in terms of almost all projections.  Indeed, let  $y=(y_1, \ldots, y_{r+1})$ be another local system of coordinates at $\point$.  Let us denote the indeterminates of \ref{eq:genprojection} in this new system of coordinates by $v^{(i)}_{\nu_1, \cdots, \nu_{r+1} }$.  Then each $v^{(i)}_{\nu_1, \cdots, \nu_{r+1} }$ is a finite linear combination of the $u^{(i)}_{\nu_1, \cdots, \nu_{r+1} }$ with coefficients determined by the coefficients of the local change of coordinates, that is of triangular form, namely  $v^{(i)}_{\nu_1, \cdots, \nu_{r+1} }$ of weight $d=\sum \nu_i$ depends only on   $u^{(i)}_{\nu_1, \cdots, \nu_{r+1} }$  of weight smaller than or equal to $d$.  
 By composing the polynomials of $\mathcal G=\{G_\mu (u) \}$ with 
$u(v)$ we obtain a finite family of polynomials $\mathcal H=\{H_\mu (v) = G_\mu (u(v))  \}$ that determines  the specializations that give Zariski equisingularity in the new system of coordinates. 
\end{rem}

\begin{rem}\label{rem:comment2} 
As Zariski shows in section 2 of \cite{Zariski1979}  it is enough to consider the generic projections of the form 
 \begin{align}\label{eq:genprojection2}
\pi_{u,i} (x) =  x_i + u^{(i)} x_{r+1} + \sum_{d\ge 2} \sum_{\nu_1+ \cdots \nu_n=d} u^{(i)}_{\nu_1, \cdots, \nu_{r+1} } x^\nu .  
  \end{align}
(This corresponds  to dividing by linear automorphisms of the target space of the projection.)
\end{rem}


\subsection{Generic linear projections}\label{ssec:glp}

One may similarly define \emph{generic linear Zariski equisingularity} as the one given by linear projections belonging to a 
Zariski  open non-empty subset of linear projections. Similarly to the case of a generic projection, it can be defined 
equivalently by the projection \eqref{eq:genprojection}, where the sum is taken only on $|\nu|=1$, with the field extension to any field containing the field of rational functions $\kk(u^{(i)}_{\nu_1, \cdots, \nu_{r+1} }, |\nu|=1)$.  
Note that, since the discriminant of a linear projection depends only of 
its kernel we may slightly simplify the picture, similarly to Remark \ref {rem:comment2},  and consider only 
the projections of the form 
 \begin{align}\label{eq:genprojectionlinear}
\pi_{u,i} (x) =  x_i + u^{(i)} x_{r+1}.  
  \end{align}

 It is not clear  whether such notion of generic linear Zariski equisingularity is preserved by non-linear local changes of coordinates,  nor whether it implies the generic Zariski equisingularity.  
Nonetheless this is the case for the dimensionality type $1$, by \cite{zariski65-S2}, and for the dimensionality type $2$, by  
Theorem \ref{thm:mainequivalence}.


\subsection{Algebraic case}\label{ssec:algebraic}
Let $f \in \kk [x_1, \ldots, x_r]$ be a reduced polynomial. Then $V = V(f)$ is a hypersurface of $\mathbf A_\kk^r$. 
Hironaka showed in \cite{hironaka79} the semi-continuity (in Zariski topology) of the dimensionality type.  It follows that the sets $V_\sigma = \{P\in V; \dimtype (V,P) \ge \sigma\}$ are Zariski closed in $V$. Hironaka's proof is based on a detailed analysis of the Weierstrass preparation theorem for $K((u))[x]$, where $K$ is a suitable base field.

Note that Hironaka's result is obvious for $\sigma=1$ because, by definition, $V_1 = \sing V$ and it is easy for $\sigma=2$.  
Indeed in the later case, we may use the generic linear projections and therefore $V_2$ is equal to the generic intersection 
of  $\pi_u^{-1}(\sing \Delta_u)$, where $u$ parameterizes the linear projections.

The precise argument is as follows.  Consider only those linear projections $\pi_u, u\in U$, that restricted to $V$ are finite. Let $X\subset V\times U$ be the Zariski closure of the union of $\pi_u^{-1}(\sing \Delta_u)$.  Denote by $X_u$ the fibres of the 
projection of $X$ on $U$.  Then $V_2 = \bigcap_u X_u$.

By Theorem \ref{thm:mainequivalence} a similar argument shows that $V_3$ is algebraic.  Indeed, it is easy to see that 
if we replace the generic Zariski equisingularity by the generic linear Zariski equisingularity then $V_\sigma$ is the generic intersection of $\pi_u^{-1}(V_{\sigma-1}(\Delta_u))$.  


\subsection{Complex analytic case}\label{ssec:analytic}
Let $M$ be a complex analytic manifold and let $V $ be a complex analytic hypersurface of $M$. 
Since in the definition of the dimensionality type it is enough to consider only the polynomial projections, the 
dimensionality type gives a canonical complex analytic stratification and the sets $V_\sigma = \{P\in V; \dimtype (V,P) 
\ge \sigma\}$ are complex analytic.  

This canonical stratification of $V$ is topologically trivial along each stratum by \cite{varchenkoizv73, varchenkoICM75} and arc-wise analytically trivial by \cite{PP17}. 
By \cite{PTpreprint}  $V$ is bi-Lipschitz trivial along the strata of codimension $1$ and by \cite{NPpreprint} or \cite{PPpreprint19} it is bi-Lipschitz trivial along the strata of codimension $2$.  


\section{Zariski equisingular families}

Zariski showed in \cite[Theorem 4.4]{zariski65-S2}  that singularities of dimensionality type $1$ are isomorphic to the total space of equisingular families of plane curve singularities. He developed the theory of such families in \cite{zariski65-S1, zariski65-S2,zariski68-S3}.  We apply a similar approach to the singularities of dimensionality type $2$ by considering 
equisingular families of surface singularities in the three-dimensional space. For this we use the notion of equisingularity of families of singularities introduced by Zariski in \cite{zariski71open} and developed by Varchenko, see 
e.g. \cite{varchenkoizv73,varchenkoizv75,varchenkoICM75}.  It is also called Zariski equisingularity, it is similarly defined by taking successive discriminants.  Thus there are two closely related notions of Zariski equisingularity: along a subspace and of a family.

In this section we recall this notion beginning with the case of plane curve singularities and then considering the general case. 
In the next section we develop the case of families of surface singularities by introducing the notion of 
\unt Zariski equisingularity.

\subsection{Zariski equisingular families of plane curve singularities}\label{subsec:ZEcurvefamilies} { \, } \\
Let $g(x,y,t)\in \,$  $\kk [[x,y,t]], t=(t_1, \ldots , t_l)$. Assume that $g(0,0,t)$ is identically zero and that $g$ is regular in the variable $y$, i.e.  $g(0,y,0) \not \equiv 0$.  We also assume that $g$ is reduced and we will study its zero set $V=V(g)=g^ {-1}(0)$ as a family of plane curve singularities $V_t\subset (\kk^2,0)$ parameterized by $t\in (\kk^l,0)$.  
By the Weierstrass preparation theorem we may assume that $g$ is of the form 
\begin{align} \label{eq:weierstrass-g}
g(x,y,t) = y^d + \sum _{i=1}^d a_i(x,t) y^{d-i}    
\end{align}
with $a_i(0,0)=0$.
\begin{defi}
We say that the family $V_t$ is 
\emph{Zariski equisingular with respect to the parameter $t$ 
and the projection $ (x,y,t)\to (x,t)$} if there is a non-negative integer $M$ such that the discriminant of $g$ satisfies 
$$D_g(x,t) = x^M \cdot unit(x,t).$$
We say that $V_t$ is 
\emph{Zariski equisingular with respect to the parameter $t$}  if there is 
such a system of local coordinates 
$(x',y',t) \in (\kk^2\times \kk^l, 0)$ (same parameter $t$) that makes 
$V_t\,$ Zariski equisingular for the projection $ (x',y',t)\to (x',t)$. 
\end{defi}

For equisingular families of plane curve singularities
 the special fiber, i.e. for $t=0$, and the generic fiber 
 are equivalent plane curve singularities, 
see  of \cite[Section 6]{zariski65-S1}. In \cite[ Section 3]{zariski65-S1}  
Zariski gave three equivalent definitions of equivalence of 
algebroid curves.  
In the complex analytic case Teissier gave 
in \cite[page 623]{teissier77} 
12 equivalent characterizations of such equivalence in families,  including the  topological equivalence.   

We say that the projection $(x,y,t) \to  (x,t)$ is  \emph{transverse} 
if its kernel is not included in the tangent cone $C_0(V)$. 
This is equivalent to the multiplicity of $V$ at $0$ being  equal 
to $d$.

\begin{rmk}\label{rk:transversality and equimultiplicity}
One important feature of the Zariski equisingularity is that it implies equimultiplicity, see \cite{zariski75} 
for general case and \cite[page 531]{zariski65-S1} for families of plane curve singularities.  
This concerns both Zariski equisingualrity along a subspace and Zariski equisingularity in families 
(we recall the latter notion in Section \ref{ssec:ZEfamilies}).  
The equimultiplicity of $V_t$, as family of reduced hypersurfaces, 
implies the normal pseudo-flatness of $V$ along $\{0\}\times \kk ^l$,  
intuitively the continuity of the tangent cone, see e.g. \cite{Hironaka69}
for the notion of normal pseudo-flatness.  Therefore, in Zariski equisingular families, 
if $V_0$ is transverse to the kernel of the projection so is $V$ and so are $V_t$ for $t$ non-zero and small.
\end{rmk}

By \cite[Theorem 4.4]{zariski65-S2}, $V$ is of dimensionality type $1$ at 
$\point$ if and only if there is a local system of coordinates $x,y,t$ at 
$\point$ such that $V$ is a Zariski equisingular family with 
respect to the parameter $t$.  Moreover, if this is the case, then 
$V$ is such an equisingular family for any 
transverse system of coordinates \cite[Theorem 7]{zariski65-S1}.

We shall use in this paper the fact that Zariski equisingular families of plane curve singularities can be parameterized by fractional power series in $x$, with parameter 
$t$, in the sense of \cite[Theorem 2.2]{PP17}. This follows from the 
Jung-Abhyankar Theorem, that is valid over any algebraically closed 
field of characteristic zero, see \cite{jung1908}, \cite{abhyankar55},  and also \cite{PR2012}.  
\begin{thm}[Puiseux with parameter]\label{thm:Jung}
Let $g\in \kk[[x,y,t]]$ be reduced and of the form \eqref{eq:weierstrass-g}.  
Then $V(g)$ is 
Zariski equisingular with respect to the projection $ (x,y,t)\to (x,t)$ if and only if there are $y_i \in \kk[[u,t]], i=1,\ldots ,d $, and strictly positive integers $n$, $k_{ij}, i<j$, such that 
\begin{align}
g(u^n,y,t) = \prod_{i=1}^d (y- y_i(u,t)) 
\end{align}
satisfying $y_i - y_j = u^{k_{ij}} unit (u,t) $.  
\end{thm}

\subsection{Zariski equisingular families.  General Case}\label{ssec:ZEfamilies}

For $x=(x_1, \ldots, x_{r+1}) \in \kk^{r+1}$ we denote 
 $x^i = (x_1, \ldots, x_i) \in \kk^i$.  

 By a \emph{local system  of pseudopolynomials in 
 $x=(x_1,... ,x_{r+1})\in \kk^{r+1}$ at $(0,0)\in \kk^{r+1}\times \kk^l$, with parameter $t\in (\kk^l,0)$}, we mean a 
 family  of power series 
  \begin{align}\label{def:pseudopolynomials}
  F_{i} (x^i,t )= x_i^{d_i}+ \sum_{j=1}^{d_i} a_{i-1,j} (x^{i-1},t)
 x_i^{d_i-j},  \quad i=0, \ldots,r+1, 
\end{align} 
with the coefficients $a_{i-1,j}$ vanishing identically for $x^{i-1}=0$. 
This includes $d_i=0$, in which case we mean $F_i\equiv 1$.    
 
\begin{defi}
Let $F\in \kk[[x_1, \ldots, x_{r+1}, t_1, \ldots , t_l]]$. 
We say that $V=V(F)$ is \emph{Zariski equisingular with respect to the parameter $t$}
 if there are $k\ge 0$ and a system of pseudopolynomials 
 $F_{i} (t, x^i )$ 
 such that 
\begin{enumerate}
\item 
$F_{r+1}$ is the Weierstrass polynomial associated to $F$.
\item
for every $i$, $k\le i\le r$,  $(F_{i+1})_{red}(t,0)\equiv 0$ and  the discriminant of $(F_{i+1})_{red}$ (or, equivalently, the first not identically equal to zero generalized discriminant of $F_{i+1}$, see Appendix IV of \cite{whitneybook}, or Appendix B of \cite{PP17}) divides 
$F_{i}$.     
\item
$F_k \equiv 1$  (and then we put $F_i \equiv 1$ 
for all $0\le i< k$).
\end{enumerate}
\end{defi}


\subsection{Families uniformly rational in a finite number 
of parameters}
\label{ssec:rational dependence}

Suppose $\tau=\{\tau_1, \ldots, \tau_s\}$ is a finite set of indeterminates, $f\in \kk(\tau) [[x]]$ and $\pi: \kk^{r+1} \to \kk^{r} $ restricted to $V (f)$ is finite.  By a \emph{specialization} 
we mean a substitution of $\tau $ by $\overline \tau \in \kk^s$ (see also the discussion in Section \ref{sec:dim-type}).  It makes sense (it is admissible) if the denominators of the coefficients 
of $f$ do not vanish at  $\overline \tau$.  In order the specialization to be well-defined (admissible) 
on a non-empty Zariski open  subset of  $\kk^s$ we consider   series $f\in \kk(\tau) [[x]]$ of a special type.  
We will omit the bar notation $\bar \tau$ whenever it is clear from the context.
Our presentation is a simplified version of the case of an infinite number of indeterminates 
developed in Sections 4 and 5 of \cite{Zariski1979}, where 
a non-empty Zariski open set  is replaced by the notion of almost all specializations.  

\begin{defn}
We say that a power series $f\in \kk(\tau) [[x]]$, 
\begin{align}
f(x,\tau) = \sum_\alpha a_\alpha(\tau ) x^\alpha,  
\end{align}
$x=(x_1, \ldots, x_{r+1}) $, 
$\tau=(\tau_1, \ldots, \tau_s)$, \emph {is uniformly rational in $\tau$}
if the  coefficients of $f$ are of the 
form $a_\alpha = b_\alpha/c^{|\alpha|}$, 
where $b_\alpha, c $ are polynomials $b_\alpha, c \in \kk[\tau]$, 
and $c\ne 0$, or, equivalently, $f\in \kk[\tau] [[\frac x c]]$. 
\end{defn}

This notion has two convenient properties that we need for the statement and the proof of Proposition \ref{prop:retionalinparameter}.  
Firstly, for a power series $f$ uniformly rational in $\tau$,  the 
 specialization $\overline \tau \to \tau$ makes sense for $\overline \tau$ from a non-empty Zariski open subset, i.e. for 
 if $c(\overline \tau) \ne 0$.  
Secondly, as we show below, it follows from a classical proof of the Weierstrass preparation theorem that the property of being
 uniformly rational in $\tau$ is preserved by taking the associated Weierstrass polynomial, and hence by taking the discriminant of $f$ as well.

\begin{lem}\label{lem:rational}
Suppose that $f\in \kk(\tau) [[x]]$ is uniformly rational in $\tau$ and is regular in $x_{r+1}$ (i.e. $f(0,x_{r+1})\not\equiv 0$).  
Then, the coefficients $a_i$  
of the Weierstrass associated polynomial of $f$ over $\kk(\tau)$
\begin{align*}
f(x,\tau) = \Big (x_{r+1}^d + \sum_{i=1}^d a_i(x', \tau) x_{r+1}^{d-i} \Big ) \ unit(x,\tau), 
\end{align*}
$x'= (x_1, \ldots , x_{r})$, are uniformly rational in 
$\tau$ (maybe with a different  denominator $c$). In particular, the discriminant 
$\D_f(x',\tau)$  of   $f$ also is  uniformly rational in $\tau$.
\end{lem}

\begin{proof}
Let us consider first  $f = \sum c_n(y) z^n\in \kk[\tau] [[y,z]]$, where $y=(y_1, \ldots , y_r)$ and $z$ is a single variable, and suppose, moreover, that there is $d\ge 1$  such that $c_i(0)=0$ for $i<d $ and $c_d(0)=1 $.
We denote by $\my_y$ the ideal of $\kk (\tau) [[y,z]]$ generated by $y_1, \ldots, y_r$ and by $\my_z$ the one generated by 
$z$.  Thus $\my := \my_y + \my_z$ is the maximal ideal of $\kk (\tau) [[y,z]]$.  
 
For a series $f$ we define $\h (f):= \sum_{n\geq 0} c_{n+d} z^n$.  
Then  $\h (f)\in 1+ \my \cap \kk [\tau ] [[y,z]]$ and therefore $\h (f)^{-1}\in 1+ \my \cap \kk [\tau ] [[y,z]]$.   Write 
$$
f= z^d \h (f) + R,  \quad \text {with} \quad R\in \my_y \cap \kk[\tau] [[y]][z]_{d-1},  
$$
that is $R$ is a polynomial in $z$ of degree at most $d-1$.  
  Define 
$f_1:= (\h (f))^{-1} f = z^d +   (\h (f))^{-1} R  $. Note that 
$(\h (f))^{-1} R \in  \my_y \cap \kk [\tau ] [[y,z]]$  and 
$$
f_1= z^d \h (f_1) + R_1,  \quad R_1\in  \my_y \cap \kk[\tau] [[y]][z]_{d-1}, \quad \h (f_1)
\in 1+ \my_y \cap \kk [\tau ] [[y,z]] .$$
The next step is similar. Define 
$f_2:= (\h (f_1))^{-1} f_1 = z^d +  (\h (f_1))^{-1} R_1 $.  Then 
$$
f_2= z^d \h (f_2) + R_1 + R_2 ,  \quad R_2\in \my_y^2 \cap \kk[\tau] [[y]][z]_{d-1}, \quad \h (f_2) \in 1+ 
\my_y^2 \cap \kk [\tau ] [[y,z]] .
$$
Recursively we set $f_{k}:= (\h (f_{k-1}))^{-1} f_{k-1} = z^d + (\h (f_{k-1}) )^{-1} (\sum_{i=1}^{k-1} R_i) ,$ 
$$
f_{k} = z^d \h (f_{k}) + \sum_{i=1}^{k} R_i, \quad 
R_{k}\in  \my_y^{k} \cap \kk[\tau] [[y]][z]_{d-1}, \quad \h (f_{k})\in 1+ 
\my_y^{k} \cap \kk [\tau ] [[y,z]] .
$$
  With $k\to \infty$ we get 
\begin{align*}
f = (\prod_{i=0}^{\infty} \h (f_i)) (z^d  + \sum_{i=1}^{\infty} R_i)  \quad \text {with} \quad \sum_{i=1}^{\infty} R_i \in 
\kk[\tau] [[y]][z]_{d-1}
\end{align*}
This shows the lemma under the assumptions $f \in \kk[\tau] [[y,z]]$ and $c_d(0) =1$.

In the general case  $f \in \kk[\tau] [[\frac y c, \frac z c]]$, $c_i(0)=0$ for $i<d $ and $c_d(0)= \frac b {c^ d} \ne 0 $.  
Set $\tilde y = \frac y {cb^{d+1}}$, $\tilde z = \frac z {cb}$.  Then $\tilde f \in \kk[\tau] [[\tilde y, \tilde z]]$, 
given by $\tilde f (\tilde y, \tilde z) := b ^ {-d-1} f(cb^{d+1}\tilde y, cb \tilde z)$ satisfies the assumptions of the case 
considered above and therefore 
\begin{align*}
\tilde f (\tilde y, \tilde z)   = \tilde u (\tilde y, \tilde z) (\tilde z^d  + \sum_{i=1}^{d} \tilde a_i (\tilde y) \tilde z^ {d-i})    \quad \text {with} \quad \tilde a_i \in 
\kk[\tau] [[\tilde y]]. 
\end{align*} 
This implies that 
\begin{align*}
 f (y, z)   & = b^{d+1}\tilde u \big (\frac y {cb^{d+1}}, \frac z {cb} \big ) \Big ( \big(\frac z {cb}\big )^d  + \sum_{i=1}^{d} \tilde a_i \big (\frac y {cb^{d+1}}\big ) \big (\frac z {cb}\big )^ {d-i} \Big)\\
 & = u(y,z) \Big (z^d  + \sum_{i=1}^{d} (cb)^i \tilde a_i \big (\frac y {cb^{d+1}}\big ) z ^ {d-i} \Big ), 
\end{align*} 
where $u(z,y) = \frac b {c^d}  \tilde u (\frac y {cb^{d+1}}, \frac z {cb}) $.  This completes the proof of lemma.    
\end{proof}

The above lemma allows us to translate, both ways, Zariski equisingularity with respect to  a global parameter 
$\tau\in U\subset \kk^s$, the theory we do not develop in this paper, to Zariski equisingularity with respect 
to a local parameter $\tau\in (\kk^s, \tau_0)$ for generic $\tau_0$.  This is given in the proposition below, 
$t$ is an additional local parameter.  

\begin{prop}\label{prop:retionalinparameter}
Suppose that $f(x,t,\tau)\in \kk(\tau) [[x,t]]$ is uniformly rational in $\tau$, 
$f(0,t)\equiv 0$, and $f$ is regular in $x_{r+1}$.  Denote by $\pi$  the projection 
$(x_1, \ldots, x_r, x_{r+1}, t)\to (x_1, \ldots, x_r, t)$. 
Then the following conditions are equivalent:
\begin{enumerate}
  \item 
$V= V(f)$, understood as an algebroid  hypersurface over $\kk(\tau)$, 
is Zariski equisingular with respect to the parameter $t \in (\kk^l,0)$ (and the projection $\pi$).  
\item
There is a non-empty Zariski open  subset of $U \subset \kk ^s$ such that 
for every $ \tau \in U$ the family 
$V_{ \tau } =\{f(x,t,\tau) =0\}$ 
is Zariski equisingular with respect to the parameter $t \in (\kk^l,0)$ (and the projection $\pi$).  
\item
There exists $ \tau_0 \in \kk ^s$ such that $V$, understood 
as an algebroid hypersurface of $(\kk^{r+1+l+s},\,  (0,0,\tau_0))$, 
is Zariski equisingular with respect to the local parameters 
$(t,\tau-  \tau_0) \in (\kk^l \times \kk^s,0)$ (and the projection $(x_1, \ldots ,x_r, x_{r+1}, t,\tau)\to (x_1, \ldots ,x_r, t, \tau)$).  
\end{enumerate}
\end{prop}

\begin{proof}
We just sketch the proof.  It follows by arguments similar to the proof 
of Theorem 5.1 of \cite{Zariski1979}, where the much more difficult case of infinitely many indeterminates 
(here corresponding to the parameter $\tau$) was considered. All three conditions can be verified by taking 
the discriminant, that is  $f$ satisfies one of these conditions 
if and only if the reduced discriminant $(\D_f)_{red}$   of  $f$ by the projection $\pi$ satisfies that condition as well.  
Here we use Lemma \ref{lem:rational} because we need that the uniform rationality is preserved by taking the discriminant.  Since, moreover, for $\tau_0$ from a non-empty Zariski open subset of $\kk^s$ the discriminant $D_f$ evaluated at $\tau_0$ equals 
the discriminant of $f(x,t, \tau_0)$ (specialization of $\tau $ to generic $\tau_0$ and taking the discriminant commutes), each implication between the three conditions can be proven by induction on $r$. 
Therefore the proof can be reduced to the case 
$r=0$.  In this case the discriminant of such $f$ is a series $D_f (t,\tau)$ uniformly rational in $\tau$ and all three conditions are equivalent to the property $D_f (0,\tau)\ne 0$ for generic $\tau$.  
\end{proof}


\section{Zariski equisingular families of surface singularities}\label{sec:main}

\subsection{Nested uniformly transverse system of coordinates}\label{ssec:nut}
Given  an algebroid hypersurface $V_0 = f_0^{-1} (0) \subset (k^3,\point )$ defined by a reduced formal power series $f_0$.  
 For $b\in \kk$ we denote by $\pi_b$ the projection 
$\kk^3\to \kk^ 2$ parallel to $(0,b,1)$, that is $\pi_b(x,y,z) = (x, y-bz)$.  
If $\pi_b$   restricted to $V$ is finite, for instance if its kernel is transverse to $V,$
then we denote by $\Delta_b$ the discriminant locus of such restriction.

\begin{defn} \label{defn:snt-generic}
We say that a local  system of coordinates $x,y,z$ at $\point$
is \emph{nested-uniformly-transverse}, or \emph{\snt} for 
short, for $V_0$  if 
\begin{enumerate}
\item
The projection $(x,y,z) \to (x,y) $ is transverse to $V_0$ at $\point$, i.e., the tangent  cone $C_\point(V_0)$  does not contain the $z$-axis. 
	\item
 The projection $(x,y) \to x $ is transverse to $\Delta_0$ at $\pi_0(\point )$.
 \item 
 The family of plane curve singularities $\Delta_b$ parameterized by 
 $b\in (\kk,0)$ is equisingular.   
\end{enumerate}   
\end{defn}

\begin{lemma}\label{lem:genericlinear==>nu}
Given a local system of coordinates at $\point,$ after performing a generic linear change of coordinates, i.e.  from a Zariski open non-empty subset of linear changes of coordinates, the new system is \snt.  
\end{lemma}

\begin{proof}
This is obvious for the transversality conditions (1) and (2).  For (3) it follows from the fact that 
any family parameterized by a finite number of parameters is generically equisingular.  
To give a more precise argument one may use (1) implies (3) of Proposition 
\ref{prop:retionalinparameter}.  Indeed, consider for simplicity only the changes of coordinates 
$(X,Y,Z) \to (X+aZ,Y+bZ,Z)$ parameterized by $(a,b) \in \kk^2$ and a variety defined over $\kk(a,b)$ by $f (X+aZ,Y+bZ,Z)=0$.  
This variety is trivially equisingular as a family, there is no parameter.  
Then it follows from (3) of Proposition 
\ref{prop:retionalinparameter} that it is equisingular if considered as a family parameterized by 
$(a,b) \in (\kk^2, (a_0,b_0))$ for  $(a_0,b_0)$ from a Zariski open non-empty subset of  $\kk^2$. 
\end{proof}

The following two theorems show that for a \unt system of coordinates the projection $\pi_0$ is generic in the sense of Zariski Equisingularity,  with respect to both linear, Theorem \ref{thm:genlin},  and nonlinear, Theorem \ref{thm:generic}, changes of coordinates. 
Since the discriminant  locus of a linear projection depends only of 
its kernel, c.f. also Remark \ref {rem:comment2}, for linear changes it suffices to show the following result. 

\begin{thm}\label{thm:genlin}
Suppose $(x,y,z)$ is a \unt system of coordinates for $V_0= V(f_0(x,y,z))$. 
 Then $V(F)$ for 
\begin{align}\label{eq:linearfamily}
F(X,Y,Z,a,b):= f_0(X+aZ,Y+bZ,Z), 
\end{align}
is Zariski equisingular with respect to the parameters $(a,b)\in (\kk^2,0)$. 
\end{thm}

Consider now the  family 
\begin{align} \label{eq:genericfamily}
F_\tau (X, Y, Z)= f_0(\impli (X,Y,Z)) ,
\end{align}
where $\impli $ is a local, not necessarily linear, analytic change of coordinates depending, moreover, analytically on a parameter $\tau\in (\kk^s,0)$, 
i.e., $\varphi \in (\kk[[x,y,z,\tau]])^3$.  

\begin{thm}\label{thm:generic}
Suppose $(x,y,z)$ is a \unt system of coordinates for $V_0= V(f_0(x,y,z))$.  Consider a local change of coordinates $\impli (X,Y,Z) $ that depends analytically on a parameter $\tau\in (\kk^s,0)$ and satisfies:  $\impli_1 (0)=\impli_2 (0)=\impli_3 (0) =0$ for all $\tau$, and $D\varphi (0)=\id$ for $\tau =0$.  
Then $V(f(\impli (X,Y,Z)))$ is Zariski equisingular with respect to 
the parameter $\tau $.  
\end{thm} 

Theorems \ref{thm:genlin} and \ref{thm:generic} are consequences of 
Theorem \ref{thm:genericwithparameter} below.


\subsection{$\nu$-transverse Zariski equisingular families of surface singularities}
\label{ssec:ZEfamiliesdim2}

Consider an analytic family of surface singularities in $\kk^3$
\begin{align}\label{eq:fwithparameter}
f(x,y,z, t)= f_t(x,y,z) \in \kk[[x,y,z,t]],
\end{align} 
where $t\in (\kk^l,0)$ is considered as a parameter.

\begin{defn} \label{defn:snt-genericZE}
We say that \emph{the family $V=f^{-1} (0)$ is \unt Zariski equisingular} (with respect to a given system of coordinates) if this  system of coordinates 
$(x,y,z) \in \kk^3$ at $\point$ is \unt for $V_0= f_0^{-1} (0)$ and the family of the discriminant loci $\Delta_{b,t}$ of the projection 
$(x,y,z,t) \to (x, y-bz,t)$ restricted to $V$ parameterized by $(b,t) \in (\kk\times \kk ^l,0)$ 
is an equisingular family of plane curve singularities.    
\end{defn}

Consider an arbitrary  analytic change of local coordinates 
\begin{align}\label{eq:change} (x,y,z,t) = \impli (X,Y,Z,T), 
\end{align}
i.e.  $\impli (X,Y,Z,T)  = (\impli_1, \impli_2, \impli_3, \impli_4)$ with $\impli_4 $ 
depending not only on $T$ but also on $X,Y,Z$.

\begin{thm} \label{thm:genericwithparameter}
Suppose that the family $V(f)$, for $f$ given by \eqref{eq:fwithparameter},
 is  \unt   Zariski equisingular.  Let \eqref{eq:change} be a local change of coordinates
that depends analytically on an additional parameter $\tau\in (\kk^s,0)$ and satisfies:  $\impli_1 (0,T)\equiv \impli_2 (0,T)\equiv \impli_3 (0,T) \equiv 0$ for all $\tau$, and $D\varphi (0)=\id$ for $\tau =0$.  

Then, the family of the zero sets of $F (X, Y, Z, T)= f(\impli (X,Y,Z,T))$ is Zariski equisingular
 with respect to  parameters  $T$ and $\tau$.  
\end{thm} 

\begin{rmk}\label{rmk:noparametertau}
For applications it is convenient to include the dependence on 
the parameter $\tau$ in the above theorem. But this is not necessary for its proof.  Indeed, the parameter $\tau$ can be added to the parameters $t$ and $T$.    
\end{rmk}

We show Theorem \ref{thm:genericwithparameter} in Section \ref{sec:proof}. 
Theorem \ref{thm:generic} follows immediately from 
Theorem \ref{thm:genericwithparameter}.  To show 
Theorem \ref{thm:genlin} it suffices to consider the local change of variables 
$(x,y,z) \to (x-az,y-bz, z)$ parameterized by $\tau =(a,b)$ and apply Theorem \ref{thm:genericwithparameter}.


\section{Proof of Theorem \ref{thm:mainequivalence}}
\label{ssec:proofMainThm}

We begin with the definition of \snt-Zariski equisingularity along $S$ at $\point$.  

\begin{defn}\label{defn:snt-genericZE-along}
 Let $(V,S, \point)$ be as in Section \ref{sec:dim-type}, that is $V=V(f) \subset (\kk^{r+1},\point )$  and $S$ is a nonsingular algebroid subspace of $\sing V$.  We suppose that $\dim S= r-2$.  We say that 
 $V$ is \emph{{}\unt Zariski equisingular along $S$ at $\point$} if there is a local  system of coordinates 
 $(x,y,z,t) \in \kk^3\times \kk^{r-2}$ at $\point$ such that $S=\{x=y=z=0\}$ and, in this system of coordinates, 
the family $V=f^{-1} (0)$ parameterized by $t$ is \snt-Zariski equisingular family in the sense of Definition \ref{defn:snt-genericZE}.  
\end{defn}

\begin{lem}\label{lem:(4)==>(1)}
If, for a system of local coordinates,  $V$ is generic linear Zariski equisingular along $S$ at $\point$ then
$V$ is \unt Zariski equisingular along $S$ at $\point$.  
\end{lem}

\begin{proof} 
This is a parameterized version of Lemma \ref{lem:genericlinear==>nu} and we may adapt the proof of the latter. 
By assumption, after a linear change of coordinates, we may suppose that $V=V(f)$, $f(x,y,z,t)$ for a local system of coordinates 
$(x,y,z,t) \in \kk^3\times \kk^{r-2}$ such that for all but finite values of $b$ the discriminant of 
$$
\pi_b(x,y,z,t) = (x, y-bz,t) 
$$ 
is Zariski equisingular as a family of plane curve singularities depending on $t$.  If, moreover, $S=\{x=y=z=0\}$ then we may argue similarly to the proof of Lemma \ref{lem:genericlinear==>nu} now using the implication (2) implies (3) of Proposition \ref{prop:retionalinparameter} to ensure that $f(X,Y+bZ, Z, t)$ is Zariski equisingular family parameterized by 
$(t,b) \in (\kk^{r-2}\times \kk, (0,b_0))$, for all but finitely many $b_0$.  

In general we may suppose that $S$ is the graph of an analytic map $(x(t),y(t),z(t))$ satisfying $x(0)=y(0)=z(0)=0$ and replace the above system of coordinates by $(x',y',z',t) = (x-x(t), y-y(t), z-z(t),t)$.   Note that 
\begin{align}
\pi_b(x',y',z',t) & = (x',y'-bz',t) = (x-x(t), y-bz -(y(t)-bz(t)), t)\\ \notag
& = \pi_b(x,y,z,t) - (x(t), y(t)-bz(t), 0), 
\end{align}
and, by assumption, the discriminant of the latter projection, for all but finite values of $b$ defines an equisingular family of plane curves singularities, parameterized by $t$, at $(x(t), y(t)-bz(t))$.  Since $S=\{x'=y'=z'=0\}$ we may apply again 
Proposition \ref{prop:retionalinparameter} to deduce that this family is Zariski equisingular with respect to 
$(t,b) \in (\kk^{r-2}\times \kk, (0,b_0))$.  

We leave to the reader to verify that for a generic linear system of coordinates the transversality conditions of Definition 
\ref{defn:snt-generic} and thus complete the proof.  
\end{proof}

Lemma \ref{lem:(4)==>(1)} shows the implication (4)$\Longrightarrow $(1) of Theorem \ref{thm:mainequivalence}.  
The equivalence (2)$\Longleftrightarrow $(3) follows form the inductive definition of the dimensionality type under generic, or equivalently almost all, projections.  It can be deduced from Proposition 5.3 \cite{Zariski1979} for instance. 
Note that both conditions mean that $S$ is the stratum of the dimensionality type stratification 
containing $\point$.  

Let us show our main result (1)$\Longrightarrow$(3).  The proof is based on Theorem \ref{thm:genericwithparameter}.  
 Suppose  that $V=f^{-1} (0)$ is \unt Zariski equisingular along $S$ and let $x,y,z,t $ be a local system of coordinates 
 as in Definition \ref{defn:snt-genericZE-along}  
 To show (3) it is enough to check Zariski equisingularity for the discriminant given by a single polynomial projection $\pi_{ u}$, with $u$ satisfying  
$ G_\mu ( u)\ne 0 $ for all $G_\mu\in \mathcal G$, as explained in Section \ref{sec:dim-type}. 
These polynomials involve only finitely many indeterminates, say $s$ of them,  that we are now going to consider as 
parameters.  Let us denote them by $\tau 
=(\tau_1, \ldots ,\tau_s)$.  Let  $ u(\tau)$ be the specialization of all indeterminates that assigns $\tau \in \kk^s$ to  these  $s$ indeterminates 
and zero to the others.  Complete $\pi_{ u(\tau)}$ to a local change of coordinates $X,Y,Z,T= h (x,y,z,t)$ so that this projection is given by $(X,Y,Z,T) \to (X,Y,T)$.  Denote the inverse of $h$ by $
\impli (X,Y,Z,T) $ as in Section \ref{ssec:ZEfamiliesdim2}. This is a local change of coordinates depending polynomially 
on $\tau$.   We may assume that $D\varphi (0)=\id$ for $\tau =0$ 
because it is enough to consider 
$\pi_{\bar u}$ of the form \eqref{eq:genprojection2}, so that we may apply Theorem \ref{thm:genericwithparameter} and infer that the family is 
Zariski equisingular with respect to the parameters $t$ and $\tau$.  To conlcude we use the implication (3) implies (2) of  Proposition \ref{prop:retionalinparameter}. This shows  (3) of Theorem \ref{thm:mainequivalence}. 

Similarly, Theorem \ref{thm:genericwithparameter} shows (1)$\Longrightarrow$(5).  
Consider two local systems of coordinates, the first one $x,y,z,t$ 
as in Definition \ref{defn:snt-genericZE-along} and an arbitrary one $X,Y,Z,T$. They are related by a local change of coordinates 
$(x,y,z,t) = \impli (X,Y,Z,T)  = (\impli_1, \impli_2, \impli_3, \impli_4)$. After a generic linear change of coordinates in 
$X,Y,Z,T$ we may assume that $D\varphi (0) = id$.  By assumption $S= \{x=y=z=0\}$.  If moreover, $S= \{X=Y=Z=0\}$ then,  by  theorem \ref{thm:genericwithparameter}, $f(\varphi(X,Y,Z,T))=0$ is a Zariski equisingular family parameterized by $T$. Theorem 
\ref{thm:genericwithparameter} allows to add local parameters, hence this is also the case for $f(\varphi(X+aZ,Y+bZ,Z,T))=0$ for 
$(a,b,T)\in (\kk^ 2 \times \kk{r-2}, 0)$.  We conclude by the implication (3) implies (2) of  Proposition \ref{prop:retionalinparameter}.    
If $S\ne \{X=Y=Z=0\}$ then one may use a similar "shift" argument as in the last part of the proof of Lemma \ref{lem:(4)==>(1)}.  
Clearly  (5)$\Longrightarrow$(4) of Theorem \ref{thm:mainequivalence}.

To complete the proof we show (3)$\Longrightarrow $(4).
Suppose that (3) holds and let $ u_0$ be such that 
$\forall_\mu G_\mu ( u_0)\ne 0 $. Let us choose a new local 
system of coordinates in which $\pi_{ u_0}$ is linear.  
Then in this system of coordinates $V$  is generic linear Zariski equisingular along $S$ at $\point$.  To see it let us denote the indeterminates of \eqref{eq:genprojection}  in this new system by $v^{(i)}_{\nu_1, \cdots, \nu_{n} }$. After Remark \ref{rem:comment1}, the finite family of polynomials $\mathcal H=\{H_\mu (v) = G_\mu (u(v))  \}$ determine  the specializations that give Zariski equisingularity in the new system of coordinates.  As a result, for $v_0 = v( u_0)$, where $v(u)$ denotes the inverse of $u(v)$, $\pi _{ v_0}$ is linear and satisfies  $\forall_\mu H_\mu ( v_0)\ne 0 $.  This implies that the set of linear projections $\pi _{ v}$ satisfying 
$\forall_\mu H_\mu ( v)\ne 0 $, that is always Zariski open, is non-empty.   
\qed

\subsection{Open questions}

\begin{question}
(Propagation of $\nu$-transversality.)\\
By equimultiplicity, see Remark \ref{rk:transversality and equimultiplicity}, for Zariski equisingular families 
the transversality of projection for the central fibre $V_0$ implies the transversality for $V$ and for the nearby fibres
$V_t$. 
We may ask whether a similar property holds for $\nu$-transversality. More precisely, given $f(x,y,z,t)$ as above such that 
the $x,y,z$ is a \unt system of coordinates for $f_0$  and that $V(f)$ is Zariski equisingular with 
respect to the parameter $t$, is  $V(f)$  \unt Zariski equisingular? 
 \end{question}

\begin{question}
Is it possible to write explicitly in terms of $f$ the polynomials of 
$\mathcal G= \{G_\mu\}$ whose non-vanishing defines the \unt projections?  
 \end{question}

 \begin{question}
 The notion of $\nu$-transversality can be generalized to higher dimensionality types.  Is there an analogue  of Theorem \ref{thm:genericwithparameter}   in higher (co)dimensions?
 \end{question}



\section{Parameterization of families contour curves}\label{sec:parameterizations}

Let $f(x,y,z, t) \in \kk[[x,y,z,t]]$, and suppose that the family  $V=f^{-1}(0)$ of surface singularities parameterized by 
$t\in (\kk^l,0)$ is \unt Zariski equisingular.  That is, the family of hypersurfaces defined by  
\begin{align}\label{eq:linearfamily-b}
F(X,Y,Z,b,t):= f(X,Y+bZ,Z,t), 
\end{align}
 is a Zariski equisingular family with respect to parameters  
 $(b,t) \in (\kk\times \kk^l,0)$.   
Let 
\begin{align}\label{contour}
\contour_F = \{F=F'_Z=0\}
\end{align}
denote the associated contour subvariety of the projection  $\pi(X,Y,Z,b,t) = (X,Y,b,t)$, that is the union of  the singular set $\Sigma_F$   of $F$ and  the closure of critical locus of $\pi$ restricted to the regular part of $V(F)$.  
Let us denote the latter subvariety   by $\polar_F$.  Geometrically, 
it is a family of polar curves parameterized by $b$ and $t$.


\begin{lem}   \label{lem:parameterization-contour}
The contour subvariety $\contour_F =\{F=F'_Z=0\}$  consists of finitely many components parameterized, 
after a ramification $X=u^n$, by 
\begin{align}\label{eq:parameterization1-t}
(u,b,t) \to (u^n, Y_i ( u,b,t), Z_i (u,b,t), b,t) , 
\end{align}
where $Y_i, Z_i \in \kk[[u,b,t]]$ and there are $k_{ij}$ such that $Y_i - Y_j = u^{k_{ij}} unit (u,b,t) $.  
\end{lem}

\begin{proof}
By Zariski equisingularity (Theorem \ref{thm:Jung} applied to $D_F$ reduced)
the components of the discriminant locus $\Delta_F$ of $F$ can be parameterized, after a ramification $X=u^n$, 
as follows $(u^n, Y_i(u,b,t), b,t)$ .  To show that 
this parameterization extends to a parameterization of the  components of 
$\contour_F $, maybe after additional ramification, we use Jung-Abhyankar Theorem, see e.g. \cite{PR2012}.  

By Zariski equisingularity two such parameterization $Y_i (u,b,t), Y_j(u,b,t) $ satisfy  
$Y_i - Y_j = u^{k_{ij}} unit (u,b,t) $.  Fix one solution $Y_{i_0} \in \kk[[u,b,t]]$  
and let $m= \max_{j\ne i_0} k_{ji_0}$.  Consider 
$$
g(u,v,z,b,t) = F (u^n,Y_{i_0} + vu^m,z,b,t) .
$$
Its discriminant equals 
$$
\D_g(u,v,b,t) = \D_F (u^n,Y_{i_0} + vu^m,b,t) =  
v^{k_{i_0}}u^{k_{i_0}m} \prod_{j\ne i_0} u^{k_jk_{ji_0}} unit(u,v,b,t)
,$$
 where $k_i$ are the multiplicities of the roots $Y_i$, and it is normal crossing in variables $u$ and $v$ so it satisfies the assumptions of the Jung-Abhyankar Theorem.  Therefore, by Proposition 2.3 of \cite{PR2012}, so does the discriminant 
 of $g(u,0,Z,b,t)=0$.  This shows that $g(u,0,Z,b,t)=0$ is an equisingular family of plane curve singularities parameterized 
 by $b,t$, that, maybe after a further ramification in the variable $u$, are of the form $Z=Z_i (u, b,t)$. This completes the proof. 
\end{proof}

\begin{cor}\label{cor:singularlocus}
$\Sigma_F = \{F=F'_Z=F'_Y=0\}$.
\end{cor}

\begin{proof}
Differentiate the composition of $F$ and the 
parameterization \eqref{eq:parameterization1-t} 
with respect to $u$, $b$ and $t$ 
\begin{align*}
0 & = nF'_X u^{n-1} + F'_Y \frac {\partial Y_i} {\partial u} 
+ F'_Z \frac {\partial Y_i} {\partial u} ,  \\
0 & = F'_Y \frac {\partial Y_i} {\partial b} 
+ F'_Z \frac {\partial Y_i} {\partial b}  + F'_{b} , \\
0 & = F'_Y \frac {\partial Y_i} {\partial t_j} 
+ F'_Z \frac {\partial Y_i} {\partial t_j}  + F'_{t_j}, \quad j=1, \ldots , l .
\end{align*}
By the second and the third equation we have $F'_{b} =F'_{t_j} =0$ on the set $\{F=F'_Z=F'_Y=0\}$.  If  moreover, $u\ne 0$ then  $F'_X =0$ by the first one. 
Thus $\{F=F'_Z=F'_Y=0, X\ne 0\} \subset \Sigma_F $.

To complete the proof we consider $X=0$.  
The discriminant $D_G$ of $G(Y,Z,b,t)=F(0,Y,Z,b,t)=0$ equals $D_F(0,Y,b,t)$.  Therefore, 
the family given by $G=0$ is an equisingular family of plane curve singularities parameterized by $b$ and $t$.  
For this family,  $\{G=G'_Z=G'_Y=0\} = \Sigma_G = \{Y=Z=0\} \subset \Sigma_F $.  
\end{proof}

We will give in Proposition \ref{prop:quadratic-t} a more precise form of the parameterization 
\eqref{eq:parameterization1-t}.  For this we need a key Lemma \ref{lem:keylemma1} that we will introduce next.   
By the formula 
\begin{align}\label{eq:F'_z}
F'_Z (X,Y,Z,b,t) = bf'_y (X,Y + bZ, Z,t) + f'_z (X,Y + bZ, Z,t), 
\end{align}
$ (X,Y,Z,b,t)\in \Sigma_F$ if and only if $
(x,y,z,t) = (X,Y+bZ,Z,t) \in \Sigma_f$, the singular set of $f$. Thus in the case when
$(u^n, Y_i, Z_i,b,t) \in \Sigma_F$ we have that 
\begin{align}\label{eq:parameterization2}
(u,b,t) \to (u^n, y_i (u,b,t), z_i (u,b,t),t), \qquad y_i  := Y_i +b Z_i, \, z_i := Z_i,
\end{align}
parameterizes a component of $\Sigma_f$.  Since the projection of $\Sigma_f$ on $\{y=z=0\}$ is finite,  both 
$y_i = Y_i + bZ_i$ and $z_i=Z_i$ are independent of $b$.

Otherwise,  when \eqref{eq:parameterization1-t} parameterizes a component of $\polar_F$, \eqref{eq:parameterization2}  parameterizes a family of polar curves in $(\kk^3,0)$  corresponding to the family of 
projections $\pi_b(x,y,z,t) = (x, y-bz,t)$ (a polar wedge in the sense of 
\cite{NPpreprint} and \cite{PPpreprint19}).  

The  below key lemma is a parameterized version of the first formula on page 278 of 
\cite{brianconhenry80} 
or of a formula on page 465 of \cite{teissierRabida82}. 

\begin{lemma}\label{lem:keylemma1}
In both cases the parameterization \eqref{eq:parameterization1-t} satisfies 
\begin{align}\label{eq:mainidentity}
Z_i =  - \frac {\partial Y_i}  {\partial b} . 
\end{align}
\end{lemma} 

\begin{proof}
We have 
\begin{align}\label{eq:polar}
F (u^n,Y_i,Z_i,b,t) = F'_Z (u^n,Y_i,Z_i,b,t) = 0 .
\end{align}
 We differentiate the first identity 
with respect to $b$ and use the second one to simplify the result 
\begin{align*}
0 =   F'_Y  \frac {\partial Y_i}  {\partial b}  + F'_Z  \frac {\partial Z_i}  {\partial b}  + 
 F'_b 
 = f'_y (u^n, Y_i+bZ_i, Z_i,t) \Bigl ( \frac {\partial Y_i}  {\partial b}  +Z_i \Bigr )
\end{align*}
If $ f'_y (u^n, Y_i+bZ_i, Z_i,t)\not \equiv 0$ then  the formula \eqref{eq:mainidentity} holds.  Note that in this case  \eqref{eq:parameterization1-t} parameterizes a component of $\polar_F$.  

If $f'_y (u^n, Y_i+bZ_i, Z_i,t)\equiv 0$ then, in addition to \eqref{eq:polar},  we have $F'_Y (u^n,Y_i,Z_i,b,t) = 0$.  Thus in this case  \eqref{eq:parameterization1-t} parameterizes a component of $\Sigma_F$ 
(see Corollary \ref{cor:singularlocus}).    
Recall that in this case  both 
$Y_i + bZ_i,$ and $Z_i$ are independent of $b$.  
This implies \eqref{eq:mainidentity}.
\end{proof}



 The induced by \eqref{eq:parameterization1-t} parameterizations of the families of polar curves and of the singular locus 
 in $(\kk^3\times \kk^l,0)$, 
$(u^n, y_i (u,b,t), z_i (u,b,t))$, are given by $z_i=Z_i $ and 
$y_i  = Y_i +b Z_i$ and satisfy, 
by \eqref{eq:mainidentity}, the following 
\begin{align}\label{eq:relations}
z_i = Z_i= - \partial Y_i/\partial b,  \quad y_i  = Y_i +b z_i, \quad \partial y_i/\partial b = b\partial z_i/\partial b.
\end{align}

 The following result appeared in the non-parameterized case in the proof of Proposition 3.4 of \cite{BNP2014}. 
 The parameterized case follows from Lemma 20.3 of \cite{NPpreprint}.  We propose below a short proof.

\begin{prop}\label{prop:quadratic-t}
There are integers $m_i \in \N$ and power series $\varphi _i (u,b,t), \psi_i (u,b, t)$ such that 
\begin{align}\label{eq:generalform-t}
& x = u^ n \\
\notag
& y_i (u,b,t) = y_i(u,0,t) + b^2 u^{m_i} \varphi_i (u,b,t) \\
\notag
& z_i (u,b,t) = z_i(u,0,t) + b u^{m_i} \psi_i (u,b,t). 
\end{align}
Moreover, either $\varphi_i \equiv\psi _i \equiv  0$ if 
\eqref {eq:parameterization1-t} parameterizes a component of $\Sigma_f$, or 
 $\varphi_i(0)\ne 0$, $\psi_i(0)\ne 0$ if it parameterizes a 
 family of polar curves.  
\end{prop}

\begin{rmk}\label{rmk:polarwedges}
The image of \eqref{eq:generalform-t} is called in \cite{BNP2014}, \cite{NPpreprint} and \cite{PPpreprint19} a polar wedge 
and it is denoted by $\PW_i$.  Geometrically it is a neighborhood of a family of polar curves  parameterized by $t$.  In this paper we say sometimes that \eqref{eq:generalform-t} is a polar wedge though it would be more correct to say that it is 
\emph{a polar wedge parameterization}. 
\end{rmk}

\begin{proof}[Proof of Proposition \ref{prop:quadratic-t}]  
The case when $(u^n, y_i, z_i,t)$ parameterizes  a component 
of $\Sigma_f$ is obvious because then $y_i, z_i$ are independent of $b$.  Therefore we suppose that  one of them, and hence by \eqref{eq:relations} both of them,   depend nontrivially on $b$.  
 Expand 
 $$\frac{\partial z_i}{\partial b} (u,b,t) = \sum_{k\ge k_0} a_k(b,t) u^k$$ 
 with $a_{k_0}(b,t)\not \equiv 0$.  To show the proposition it suffices to show that   
  $a_{k_0} (0,0)\ne 0$.

Suppose,  by contradiction, that $a_{k_0} (0,0)=0$.  Then there is a  solution 
$(b(u), t(u))$, with $(b(0), t(0))=0$, of the equation $\frac{\partial z_i}{\partial b} (u,b, t)=0$.  After replacing $u$ by $u^n$, if necessary,  we may suppose that this solution is analytic.  
By the last identity of \eqref{eq:relations}, $(b(u),t(u))$ 
also solves $\frac{\partial y_i}{\partial b} =0$. Recall that $f'_z + bf'_y$  vanishes identically on the polar wedge \eqref{eq:parameterization1-t}. Hence $\frac{\partial }{\partial b} (f'_z + bf'_y)$ also vanishes on \eqref{eq:parameterization1-t} and for 
$(u,b,t) = (u,b(u), t(u))$ we get  
\begin{align*}
0= \frac{\partial }{\partial b} (f'_z + bf'_y)  = 
(f''_{yz}  + bf'_{yy}) \frac{\partial y}{\partial b}  + (f''_{zz}  
 + bf'_{yz} )\frac{\partial z}{\partial b} + f'_y= f'_y .
\end{align*}

Therefore, by Corollary \ref{cor:singularlocus}, 
\eqref{eq:parameterization1-t} parameterizes in this case a component of $\Sigma_f$.  
\end{proof}

We conclude by two lemmas.  In the first one we compare two branches $Y_i(u,b)$, $Y_j(u,b)$.  
Recall that by Theorem \ref{thm:Jung} there are  $k_{ij}\in \N$ 
 such that 
$Y_i(u,b,t)- Y_j(u,b,t) = u^{k_{ij}} unit (u,b,t)$.  By \eqref{eq:relations} this implies the following result.

\begin{lemma}\label{lem:contacts2}
\begin{align}\label{eq:generalformpairs-t}
& y_i(u,b,t)- y_j(u,b,t) = u^{k_{ij}} unit (u,b,t), \quad z_i(u,b,t)- z_j(u,b,t) = O(u^{k_{ij}})
\end{align}
\end{lemma} 

\begin{lemma}\label{lem:transversality}
For every $i$
\begin{align}\label{eq:tranversality}
y_i(u,b,t) = O(u^n), \quad z_i(u,b,t) = O(u^n).
\end{align} 
Hence always $m_i\ge n$.  We always have  $k_{ij} \ge n$ and if $m_i \ne m_j$ then $k_{ij} \le \min \{m_i, m_j\}$.   
\end{lemma}

\begin{proof}
By condition (2) of Definition \ref{defn:snt-generic}, $Y_i(0,u,0) = y(0,u,0) = O(u^n)$.  By Remark 
\ref{rk:transversality and equimultiplicity}  transversality is preserved in Zariski equisingular families. Hence $Y_i(u,b, t) = O(u^n)$.  
Now \eqref{eq:tranversality} follows from \eqref{eq:relations}.  The remaining claims are easy.  
\end{proof}

{}

\section{Proof of Theorem \ref{thm:genericwithparameter} }\label{sec:proof}

We assume  $f(x,y,z,t) \in \kk[[x,y,z,t_1, \ldots, t_l ]]$ satisfies the assumptions of Theorem \ref{thm:genericwithparameter}. 
Thus we suppose $f$ reduced, $f(0,t)\equiv 0$, and that $V(f)$ is \unt Zariski equisingular, see Definition \ref{defn:snt-genericZE}.  We may assume, moreover, that $f$ is in the Weierstrass  form 
\begin{align} \label{eq:weierstrass-f}
f(x,y,z,t) = z^d + \sum _{i=1}^d a_i(x,y,t) z^{d-i}.   
\end{align}

Let $ (x,y,z,t)= \varphi (X,Y,Z,T)$ be a local change of coordinates as in the assumptions:  $\impli_1 (0,T)\equiv \impli_2 (0,T)\equiv \impli_3 (0,T) \equiv 0$, and $D\varphi (0)=\id$.  
We show that the family of the zero sets of 
$$F (X, Y, Z, T): = f(\impli (X,Y,Z,T))$$
 is Zariski equisingular
 with respect to  parameter  $T$ (by Remark \ref{rmk:noparametertau} it is not necessary to consider the additional parameter $\tau$).  

First note that Zariski equisingularity depends only on the projection 
$(x,y,z,t) \to (X,Y,T)$, given by the corresponding coordinates of 
$\varphi ^ {-1} $.  Because the discriminant of this projection restricted to $f^{-1}(0)$ is independent on the $Z$-coordinate of $\varphi ^ {-1} $ we may suppose that $Z=z$, or equivalently $\impli_3 (X,Y,Z,T) =z$.     
Secondly, Zariski equisingularity of a family is preserved by regular  
isomorphisms of the parameter space.  Thus, because  $\varphi$ preserves the $0\times T$-subspace, after a reparameterization we may suppose that 
$\impli (0,T) = (0,T)$.  

Let $\pi(X,Y,Z,T) =  (X,Y,T)$. 
The contour set $\contour_F $ of $F^{-1} (0)$ by $\pi$
 is given by the equations :
\begin{align}\label{eq:equationforpolar}
F=0 \quad \text{ and } \quad F'_Z = f'_z\circ \impli + \impli'_{1,Z} \cdot (f'_x \circ \impli) + 
\impli'_{2,Z} \cdot (f'_y \circ \impli) + \impli'_{4,Z}\cdot (f'_t \circ \impli) =0 .
\end{align}
Its image $\pi (\contour_F )$ is the discriminant locus 
$\Delta_F $. To prove Theorem \ref{thm:genericwithparameter} one has to show that $\Delta_F $ is an equisingular family of curves parameterized by $T$.  For this we show that $\Delta_F $ satisfies the criterion given by Theorem  \ref{thm:Jung}, that is after a ramification 
 $\Delta_F $ it is the union of finitely many analytic branches with the contact orders independent of $T$.

We proceed in two steps. 
First we compute the image  $\impli (\contour_F)$, that is the contour set 
of $f^{-1} (0)$ of the projection $(x,y,z,t) \to (X,Y,T)$, more precisely of its pull-back by the polar wedges 
parameterizations \eqref{eq:generalform-t}.  (For uniformity of presentation, we consider at the same time the components of the singular set $\Sigma _f=\varphi(\Sigma_F)$, parameterized by \eqref{eq:generalform-t} with $\varphi_i\equiv\psi_i \equiv 0$.)  
We show that each polar wedge $\PW_i$ contains a component of $\impli (\contour_F)$ and 
that the contact orders  between these components, 
understood as family of curves parameterized by $t$, are independent of 
$t$ and are equal to $k_{ij}/n$ in the notation of \eqref{eq:generalformpairs-t}.

In the second step we use the inverse map $\impli^{-1}$ to parameterize 
the corresponding components of  $\contour_F = \impli^{-1}(\impli (\contour_F))$.  
Using the assumptions on $\varphi$ and the IFT we show that the 
contact orders between their projections, i.e. 
 the branches of $\Delta_F $, are independent of $T$.  
Finally we will note that these computations also imply that all the components of $\impli (\contour_F)$ are included in the union of the polar wedges $\PW_i$ and the singular locus. This  will complete the proof.

We begin by computing  $\impli (\contour_F)$. By \eqref{eq:equationforpolar} it is given by the system of equations 
\begin{align}\label{eq:relations2}
f=0 \quad \text{ and } \quad f'_z + \derimpli_1 f'_x  + 
\derimpli_2 f'_y  + \derimpli_4 f'_t =0 ,
\end{align}
where we denote $\derimpli_i = \impli'_{i,Z}\circ \impli^{-1}$.  
The $\derimpli_i$ are power series in $x,y,z,t$.  
By the assumption $D\varphi (0)=\id$,  we have    
$\derimpli_1(0)=\derimpli_2(0)=\derimpli_4(0)=0$.  

Because $f$ is \unt Zariski equisingular, each component the contour set of the projection $(x,y,z,t) \to (x,y,t)$ is included in a polar wedge (or it is a component of $\Sigma_F$). On polar wedges, that is after composing with \eqref{eq:generalform-t}, $f'_x, f'_y , f'_z , f'_t$ are related by the following formulas 
\begin{align*}
& 0=\partial x/\partial u f'_x + \partial y/\partial u f'_y + \partial z/\partial u f'_z =   nu^{n-1} f'_x + 
(\partial y/\partial u - b \partial z/\partial u ) f'_y\\
& 0=\partial x/\partial t f'_x + \partial y/\partial t f'_y + \partial z/\partial t f'_z +f'_t =  
(\partial y/\partial t - b \partial z/\partial t ) f'_y + f'_t \\
& 0 =f'_z+bf'_y .
\end{align*}

Fix a polar wedge $\PW_i$ parameterized by \eqref{eq:generalform-t}.  Using the above formulas we may express on $\PW_i \cap \impli (\contour_F)$, $b$ as a function of $u,t$.  For this we replace  first $f'_x, f'_z, f'_t$ in \eqref{eq:relations2} 
\begin{align*}
0 = - (1/n) u^{1-n}  (\partial y/\partial u - b \partial z/\partial u ) \derimpli_{1} f'_y + 
(\derimpli_{2}- b - \derimpli_{4}(\partial y/\partial t - b \partial z/\partial t )) f'_y  , 
\end{align*}
obtaining therefore an implicit equation on $b(u,t)$   
\begin{align}\label{eq:impli-b} 
b 
= \derimpli_{2} + (b \partial z/\partial t -\partial y/\partial t  )
\cdot \derimpli_{4}  + (1/n) u^{1-n}  (\partial y/\partial u - 
b \partial z/\partial u ) \cdot \derimpli_{1} ,
\end{align}
where, on the polar wedge $\PW_i$, by  
$\derimpli_k (u,b,t)$ we mean $\derimpli_k (u^n, y_i (u,b,t), z_i (u,b,t),t)$,  
and therefore 
\begin{align}
\frac{\partial \derimpli_k}{\partial b} = O(u^n).  
\end{align}
We also have 
\begin{align}
\frac {\partial } {\partial b} (u^{1-n} (\partial y/\partial u - 
b \partial z/\partial u ) ) = O(1).
\end{align}
Since, moreover, $\derimpli_1 (0)=0$, the $\partial / \partial b$ of the RHS of \eqref{eq:impli-b} vanishes on $(u,t)= 0$. By the  IFT,  
there is a unique solution of \eqref{eq:impli-b},    
 $b\in \kk[[u,t]]$ vanishing at $(u,t)= 0$.    
We denote it by $b_i (u,t)$.  Thus each component of $\impli (\contour_F)$ is  parameterized by 
\begin{align}\label{eq:generalform-phi}
& x = u^ n \\
\notag
& y_i(u,t) = y(u,b_i (u,t), t)\\
\notag
& z_i(u,t) = z(u,b_i (u,t) , t)
\end{align}

If we compare two such functions $b_i$ and $b_j$, 
for two different polar wedges, their difference is divisible by $u^ {k_{ij}}$.  
Let us show the latter claim.  First note that by Lemma \ref{lem:contacts2}
\begin{align}\label{eq:generalform-phi-diff}
& y_i(u,t)-y_j(u,t)=  y_i(u,b_i (u,t), t)-y_j(u,b_j (u,t), t) = O(u^n(b_i-b_j)) + O(u^{k_{ij}})
\\
\notag
& z_i(u,t)-z_j(u,t)=   O(u^n(b_i-b_j)) + O(u^{k_{ij}}). 
\end{align}
Therefore, for $k=1,2,4$, 
\begin{align*}
\derimpli_{k} (u^n,y_i(u,t), z_i(u,t),t) 
- \derimpli_{k} (u^n,y_j(u,t), z_j(u,t),t)  
 = O(u^n(b_i-b_j)) + O(u^{k_{ij}}). 
\end{align*}
Finally by taking the difference of two formulas 
\eqref{eq:impli-b}, the ones on $\PW_i$ and on $\PW_j$, with $b=b_i$ for the first one and $b=b_j$ for the second we get 
\begin{align*}
b_i(u,t)-b_j(u,t) 
 = o(b_i-b_j) + O(u^{k_{ij}}),  
\end{align*}
that implies $b_i(u,t)-b_j(u,t) = O(u^{k_{ij}})$.  
This, together with \eqref{eq:generalform-phi-diff}, show that we have an analog of Lemma \ref{lem:contacts2}
\begin{align}\label{eq:newgeneralformpairs}
& y_i (u,t)- y_j (u,t) = u^{k_{ij}} unit  \\
\notag
& z_i(u,t)- z_j(u,t) = O(u^{k_{ij}}) ,
\end{align}
and this shows the constancy (independence of $t$) of contact orders for $\impli (\contour_F)$.  

Next we show that the formulas \eqref{eq:generalform-phi} and \eqref{eq:newgeneralformpairs} for $\impli (\contour_F)$ induce similar formulas for 
$\contour_F$.  Denote $\impli ^{-1} =(h_1,h_2, h_3, h_4)$.  Then $Z= h_3(x,y,z,t)=z$.   On the components of $\impli (\contour_F)$, using the fact that 
$D\impli ^{-1}(0)=\id $, and $h_i (0,0, 0, t)=0$ for $i=1,2,4$ we get 
\begin{align*}
& X= h_1 (u^n,y_i (u,t), z_i (u,t), t) = u^n (1+ o(1)) \\
& T= h_4 (u^n,y_i (u,t), z_i (u,t), t) = t(1+ o(1))  .
\end{align*}
Let $X=U^ n$.  For any $n$-th root of unity we may define 
$U\in \kk[u,t]$ and this series satisfies $U := \theta u (1+ o(1))$ (one may choose for instance $\theta=1$).  Then by the Inverse Function Theorem 
$(U,T)$  is an invertible mappings of $(u,t)$.  Therefore we may reparameterize the components of $\contour_F$
as 
$$(X,Y,Z,T)= (
U^ n , Y_i(U,T), Z_i(U,T),T)$$
that satisfies  formulas analogous to \eqref{eq:generalform-phi} and \eqref{eq:newgeneralformpairs}.  
In particular it shows that $Y_i(U,T) -Y_j(U,T) 
=U^{k_{ij}}unit$ as claimed.  

To complete the proof we note that all the components of $\impli (\contour_F)$ are included in the union of the polar wedges 
$\PW_i$. Indeed, the Weierstrass polynomial associated to the discriminant $\D_F$, as a polynomial in $Y$, is of the same degree as the one associated to $\D_f$ with respect to the variable $y$.  This follows from the transversality conditions  (1) and (2) of Definition \ref{defn:snt-generic} that are preserved by small perturbations.  Therefore all the solutions of $D_F=0$ 
are among the functions $Y_i(U,T)$ considered above. This completes the proof.  

\bibliographystyle{siam}
\bibliography{ZE}

\end{document}